\documentclass{amsart}
\usepackage{graphicx}

\newcommand\RR{\mathbb{R}}

\newcommand{\FS}{{\mathcal{ F}_S}}

\newcommand{\FJ}{{\mathcal{ F}_J}}
\newcommand{\PL}{{\mathcal{ P}_{\Lambda}}}
\newcommand{\PS}{{\mathcal{ P}_{S}}}
\newcommand{\OS}{{\mathcal{ O}_{S}}}
\newcommand{\diag}{\hbox{\rm diag}\,}
\newcommand\GS[1]{{\big[ {#1} \big]_Q}}
\newcommand\GST[1]{{\big[ {#1} \big]_Q^T}}
\newcommand\TS[1]{{\big[ {#1} \big]_R}}
\newcommand\TSI[1]{{\big[ {#1} \big]_R^{-1}}}

\begin{document}

\title{Slices of matrices --- a  scenario for spectral theory}

\author{Ricardo S. Leite}
\address{Departamento de Matem\'atica, PUC-Rio, R. Mq. de S. Vicente 225, 22453-900, Rio de 
Janeiro Brazil}
\email{rsl@mat.puc-rio.br}
\thanks{The authors were supported by CNPq, FINEP and FAPERJ.}

\author{Carlos Tomei}
\email{tomei@mat.puc-rio.br}


\date{February 5, 2002}


\begin{abstract}
Given a real, symmetric matrix $S$, we define the slice $\FS$ through $S$
as being the connected component containing $S$ of two orbits under
conjugation: the first by the orthogonal group, and the second by
the upper triangular group. We describe some classical constructions in 
eigenvalue computations and integrable systems which keep slices invariant ---
their properties are clarified by the concept. We also parametrize the
closure of a slice in terms of a convex polytope.  
\end{abstract}

\maketitle

\section*{The basic definition}

Let $S$ be a real $n \times n$ symmetric matrix with simple spectrum
$\sigma(S) = \{\lambda_1 > \ldots > \lambda_n\}.$ When is a matrix 
simultaneously an orthogonal and an upper triangular conjugation of $S$?
More precisely, we consider the \emph{slice}
$\FS$ through $S$, defined to be the intersection of the sets
$\{Q^T S Q$,  for arbitrary real orthogonal matrices $Q$  with $\det Q = 1 \}$ 
and $\{R S R^{-1}, $
for arbitrary real upper triangular matrices $R$ with $\diag R > 0\}.$

How large is a slice? Clearly, if $Q^T S Q = R S R^{-1}$ then $ SQR = QRS$, and hence,
since $S$ has simple spectrum, $QR = f(S),$ for some (real) polynomial $f$.
The matrices $Q$ and $R$ are uniquely determined from $f(S)$ if
$f(S)$ is invertible: this is the standard \emph{QR-factorization} of a matrix (\cite{G}).
For convenience, we write $f(S) = QR = \GS{f(S)} \TS{f(S)}$.

Slices have been appearing in disguised form in the literature of 
numerical analysis and integrable systems, and questions about the geometry of slices 
have come up intermittently. By putting these three aspects side by side, we expect
to convince the reader that the concept is indeed a natural one.

\section*{Slices and the computation of eigenvalues}

Francis, in his fundamental work on the $QR$ algorithm
(\cite{Fr}) considered the following map between matrices. Take $S$ invertible symmetric,  
factor $S = Q R$ and define the\emph{ $QR$ step} $S' = R Q$.
It is clear that $S'$ is symmetric with the same spectrum of $S$, since $S' = Q^T S Q$. 
But more is true: $S' =R S R^{-1}$. In particular, as Francis had already pointed out,
$S$  and $S'$  have the same bandwidth (the reader should compare this argument with
the usual one \cite{P}). Also, $S'$ belongs to the slice through $S$:
in the notation for the elements of $\FS$ presented in the introduction, 
$S'$ is associated to the function $f(x) = x$.

Numerical analysts also know well that $S^{(k)} = \GST{ S^k} S \GS{ S^k}$ equals the matrix
obtained by applying $k$ times the $QR$ step starting from $S$: in 
other words, $S^{(k)}$ is the $k$-th step of the $QR$ iteration starting from $S = S^{(0)}$. Again,
$S^{(k)}$ belongs to $\FS$, and is associated to $f(x) = x^k$. 

Instead, by taking $f(x) = x^{1/k}$, the resulting 
matrix $S^{(1/k)}$ is, in  a precise sense, the $1/k$-th $QR$ step. 
Taking the limit (\cite{T1})
$$ \lim_{k \mapsto \infty} (S^{(1/k)} - S^{(0)})k= 
\lim_{\epsilon\mapsto 0} \frac{S^{(\epsilon)} - S}{\epsilon}$$
yields a vector field $\dot{S} = X(S)$,
whose solution $S(t)$ starting at $S(0) = S$ at time $k$ equals $S^{(k)}$! 
The computation of this limit appears in a number of arguments 
in the subject. The expression for $S^{(\epsilon)}$ is
$$S^{(\epsilon)} = \GST{S^{\epsilon}} S \GS{S^{\epsilon}}, \hbox{ where  }
S^{\epsilon} = \GS{S^{\epsilon}} \TS{S^{\epsilon}}.$$ 
Evaluating the derivatives in $\epsilon$, we learn that
$$\frac{ d\ }{d\epsilon} S^{(\epsilon)} =
 (\frac{d\ }{d\epsilon}\GST{S^{\epsilon}}) S\GS{S^{\epsilon}} +
\GST{S^{\epsilon}} S (\frac{d\ }{d\epsilon}\GS{S^{\epsilon}}),$$
$$(\log S) S^{\epsilon} = 
(\frac{d\ }{d\epsilon}\GS{S^{\epsilon}}) \TS{S^{\epsilon}} +
\GS{S^{\epsilon}}(\frac{d\ }{d\epsilon}\TS{S^{\epsilon}}).$$
From the last equation,
$$ \GST{S^{\epsilon}} (\log S) S^{\epsilon}\TSI{S^{\epsilon}} = 
\GST{S^{\epsilon}}(\frac{d\ }{d\epsilon}\GS{S^{\epsilon}})  +
(\frac{d\ }{d\epsilon}\TS{S^{\epsilon}})\TSI{S^{\epsilon}},$$
which obtains
$$ \GST{S^{\epsilon}} (\log S) \GS{S^{\epsilon}} \TS{S^{\epsilon}} \TSI{S^{\epsilon}}  
=\log S^{(\epsilon)} = 
\GST{S^{\epsilon}}(\frac{d\ }{d\epsilon}\GS{S^{\epsilon}})  +
(\frac{d\ }{d\epsilon}\TS{S^{\epsilon}})\TSI{S^{\epsilon}}.$$
Now evaluate the derivatives at $\epsilon = 0$: we must have 
$I = S^0 = I . I = \GS{S^0} \TS{S^0}$ and $S^{(0)} = S$. The equation above yields
$$ \log S = \frac{d\ }{d\epsilon}\GS{S^{\epsilon}}|_{\epsilon=0} + 
\frac{d\ }{d\epsilon}\TS{S^{\epsilon}}|_{\epsilon=0}.$$ 
The two terms in the right hand side of the last equation are 
special matrices: they are respectively skew symmetric and upper triangular.
Consider the (unique, linear)
decomposition of a matrix $M = \Pi_a M + \Pi_u M$ as a sum of a skew symmetric and an
upper triangular matrix. Then, from the expression for the derivative of $S^{(\epsilon)}$,
$\frac{d\ }{d\epsilon}{\GS{S^{\epsilon}}}|_{{\epsilon = 0}} = \Pi_a \log S$. 
The vector field which interpolates
the $QR$ iteration then is $ X(S)= [ S , \Pi_a \log S]$. 

The iterations and flows defined above lie in $\FS$. 
In particular, both preserve the eigenvalues of the initial condition, its symmetry and 
its bandwidth.

For numerical analysts, the
asymptotic behavior of the $QR$ iteration, well known to Francis, is of capital importance.
Say, for example, that $J$ is an arbitrary Jacobi matrix 
(i.e., a real, tridiagonal matrix whose  entries 
$J_{k,k+1} = J_{k+1,k},k=1,\ldots,n-1$ are strictly positive). 
Starting with $J$, the $QR$ iteration converges to
a diagonal matrix $D$. Not only $D$ and $J$ have the same spectrum but $D$ must
lie in the closure of $\FJ$, the slice through $J$. Steps of $QR$ type related to
different functional parameters $f$ give rise to iterations which, starting from $J$,
always converge to diagonal matrices, with diagonal entries consisting of the (distinct)
eigenvalues of $J$ in an arbitrary order. The closure of the slice $\FJ$  clearly ought to 
contain additional points. In a nutshell, they correspond to \emph{reducible} Jacobi matrices,
i.e., matrices for which some entry $J_{k,k+1} = J_{k+1,k}$ is zero. This will be 
explained in the sequel.

\section*{Slices and the Toda flows}

We now describe briefly the Toda lattice, which is an integrable system for which slices
appear as natural phase spaces. Consider $n$ particles on the line with positions
$x_k$ and velocities $y_k$ evolving under the Hamiltonian

$$ H = \frac{1}{2} \sum_{k=1}^n y_k^2 + \sum_{k=1}^{n-1} e^{(x_k - x_{k+1})}.$$ 

This Hamiltonian was introduced as a model of wave propagation in one dimensional
crystals (\cite{To}).
It was Flaschka's remarkable discovery (\cite{F}) that this dynamical system is
equivalent to the matrix differential equation
$$ \dot{J} = [ J , \Pi_a J ],$$
where $J$ is the Jacobi matrix with nonzero entries 
$$J_{k,k} = - y_k / 2 \ \hbox{    and     }\  J_{k,k+1} = J_{k+1,k} = e^{(x_k - x_{k+1})/2}/2.$$

The fact that the differential equation is in the so called
\emph{Lax pair form} (\cite{L})  implies that
the evolution $J(t)$ is actually an orthogonal conjugation of the initial condition
$J(0)$. But again more is true: the evolution stays within the set of Jacobi matrices,
as we should expect, since there are no other physically significant variables in
the problem (velocities essentially are the diagonal entries, and distances between
particles with neighboring indices correspond to the off-diagonal entries). 
Indeed, it is not hard to check that
$$J(t) = \GST{e^{t J(0)}} J(0) \GS{e^{t J(0)}} 
= \TS{e^{t J(0)}} J(0) \TSI{e^{t J(0)}},$$ the celebrated solution of the Toda lattice
by factorization (\cite{S1}). To check this formula, proceed as in the computation of
the limit of the previous section. Thus, again, the orbit $J(t)$ lies within the 
slice through $J(0)$. 

Moser also computed the asymptotic behavior of the Toda lattice (\cite{M}): for $f(x) = x$,
the orbit $J(t)$ starting from a Jacobi matrix $J(0)$ converges to a diagonal matrix,
with eigenvalues disposed in decreasing order. This is in accordance with the
relationship between Toda and $QR$ and has a natural physical interpretation. 
Diagonal entries are velocities: in the long run, particles move apart and tend to
undergo uniform motion, each with speed given by a different eigenvalue. Faster
particles move ahead, explaining the orderd outcome of the eigenvalues along the
diagonal found in applications of the $QR$ iteration. The 
remarkable fact that asymptotic speeds both at $- \infty$ and $+ \infty$ are the same
will not be relevant to us: this is strong evidence of the integrability of the system ---
each orbit 'remembers' this data. Other matrices belonging to the boundary of the
slice through $J(0)$ correspond to asymptotic behavior in which the system of particles
breaks into essentially disconnected components, the so called \emph{clustering}.

For a symmetric matrix $S$ and an arbitrary function $g$, the matrix equation 
$ \dot{S} = [ S , \Pi_a g(S)]$ admits a similar
solution by factorization: conjugate the initial condition by $\GS{e^{t g(S(0))}}$.
It is this last fact which historically was responsible for relating
the $QR$ iteration and Toda flows (\cite{S2}, \cite{DNT}):  
when $g(x) = \log (x)$, the resulting differential equation gives rise to the orbits 
interpolating the $QR$ iteration. For applications of these differential equations 
to numerical analysis, the reader may consult \cite{DLT}.

\section*{Geometric aspects of slices}

Detailed study of variables for which the Toda flow becomes especially simple led
Moser (\cite{M}) to a parametrization of Jacobi matrices. It turns out that Jacobi
matrices have simple spectrum and its eigenvectors always have nonzero first coordinates ---
in particular, they can be normalized so as to have (strictly) positive first coordinates.
Moser proved that the map taking Jacobi matrices to $n$-uples of distinct real numbers
(the eigenvalues) and a point in the first octant of the unit sphere in $\RR^{n}$
(the first coordinates of the normalized eigenvectors) is a diffeomorphism. He then
showed that an appropriate choice of Toda flow (i.e., of functional parameter $g$) 
gives rise to orbits joining any two Jacobi matrices with the same spectrum.
Thus, given a Jacobi matrix $J$,
$\FJ$ is the set of all Jacobi matrices with the same spectrum as $J$.
(We remind the reader that Jacobi matrices always have simple spectrum.) From
Moser's parametrization, \emph{Jacobi slices} are diffeomorphic to $\RR^{n-1}$. 

More generally, getting back to the definition of slices, we have seen that a matrix 
$S' \in \FS$ is of the form $S' = Q^T S Q = R S R^{-1}$, where $f(S) = QR$
for some function $f$. Keeping account of the requests on $Q$ and $R$, one may
obtain a simple coordinatization of a slice. 
A matrix $S$ is \emph{irreducible} if $S$ has no invariant subspace 
generated by a subset of the canonical vectors $e_1, \ldots, e_n$. Indeed, there
is a diffeomorphism between positive polynomials up to scalar multiplication and elements
of $\FS$, for an irreducible symmetric matrix (\cite{LT}). 

Which matrices are in the boundary $\partial \FJ$ of the slice through a 
Jacobi matrix $J$? Numerical analysts knew that the diagonal
matrix $\lambda = (\lambda_1,\ldots,\lambda_n)$ and its permuted counterparts
$\lambda_\pi = (\lambda_{\pi(1)},\ldots,\lambda_{\pi(n)})$. 
As we shall see, the simplest possible way of combining such 
$n!$ points to form a reasonable boundary turns out to be the 
topological description of $\partial \FJ$.
Let $\PL$ be the \emph{permutohedron},
obtained by taking the convex closure of the points of the form $\lambda_{\pi} \in \RR^n$.
Notice that actually $\PL$ has an interior of dimension $n-1$, since the sum of the coordinates
of any point in $\PL$ equals the trace of $J$. 
The permutohedron $\PL$ is homeomorphic to the closure
of $\FJ$, as was proved in \cite{T2} with combinatorial arguments. 

The figure below describes the situation for a $3 \times 3$ Jacobi matrix $J$ with 
spectrum $\sigma(J) =\{4,2,1\}$. On the left is a topological representation of the
closure of the slice $\FJ$. Interior points of the hexagon are points in $\FJ$. 
The boundary consists of two kinds of points: the six diagonal matrices, corresponding
to the vertices, and the matrices which form the edges, which have either entry
$(1,2)$ or $(2,3)$ equal to zero. On the right, the permutohedron associated to 
$\sigma(L)$ is projected on the $(x,y)$ plane. Notice that a map taking vertices
to vertices, say, taking $\diag(1,2,4)$ to $(1,2,4)$ and $\diag(2,1,4)$ to $(2,1,4)$
may not behave so na\"\i vely on all vertices, for continuity reasons.

\begin{figure}[tb]
\includegraphics*[width=12cm]{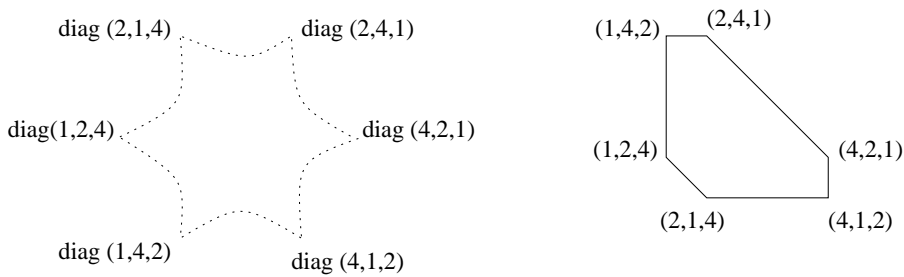}
\caption{Slice and permutohedron for $\sigma(J)=\{4,2,1\}$}
\label{fig1}
\end{figure}

The permutohedron associated to a Jacobi $4 \times 4$ matrices is also drawn below.
Faces have a clear meaning: there must be four hexagons associated
to matrices for which entry $(1,2)$ equals zero (and hence contained a Jacobi
$3 \times 3$ block of fixed spectrum), six quadrilaterals corresponding to matrices 
with entry $(2,3)$ equal to zero, and four more hexagons for matrices with 
entry $(3,4)$ equal to zero. There are also $4!$ vertices.

\begin{figure}[tb]
\includegraphics*[width=8cm]{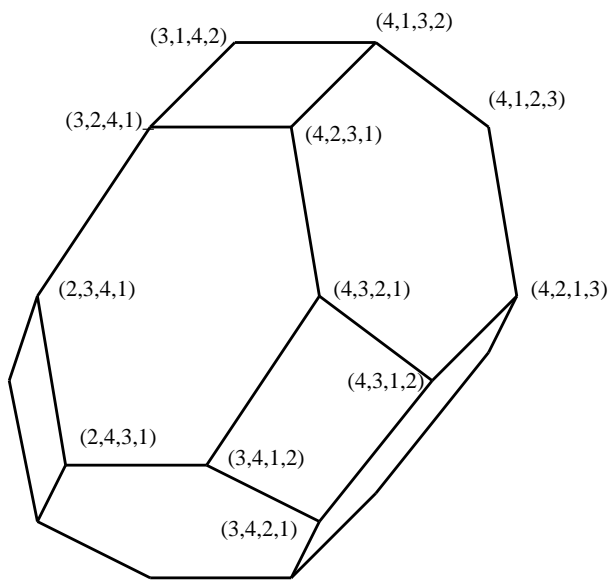}
\caption{The permutohedron for $\sigma(J)=\{4,3,2,1\}$}
\label{permu1}
\end{figure}

The permutohedron is familiar to spectral theorists, from the well known Schur-Horn
theorem (\cite{MO}). Let $S$ be a symmetric matrix with simple spectrum 
$\sigma(S) = \{\lambda_1,\ldots,\lambda_n\}$,
and consider the set $\OS = \{ Q^T S Q, Q \hbox { orthogonal } \}$. Then the map
$S \mapsto \diag S$ is surjective from $\OS$ to $\PL$. Such a map is highly not
injective. Slices are in a sense minimal sets  in $\OS$ in bijection with permutohedra. 

The presence of a convex polytope raised the possibility of relating slices to moment maps
of symplectic toroidal actions (\cite{A},\cite{GS}). Indeed, following this route, Bloch
Flaschka and Ratiu obtained an \emph{explicit} diffeomorphism between the closure of
$\FJ$ and $\PL$, which we now describe. Let $S$ be a symmetric matrix with simple spectrum, 
and consider the spectral decomposition 
$S = Q^T \Lambda Q$, where, as usual, $Q$ is orthogonal and $\Lambda$ is diagonal.
Now, order the diagonal entries of $\Lambda$ in descending order: 
$Q$ is then defined up to a choice of sign for each 
column of $Q^T$ (since they are normalized eigenvectors of $S$). Thus, 
$S = Q^T D \Lambda D Q$, for some diagonal matrix $D$ of signs. 
The matrix $ D Q \Lambda Q^T D$ is dependent on $D$, but
its diagonal is not! The map $S \mapsto \diag( Q \Lambda Q^T)$ is the \emph{BFR map}. 
Its restriction to the closure of a Jacobi slice $\FJ$ is the required diffeomorphism to
$\PL$ (\cite{BFR}).

The proof of the result above makes use of sophisticated machinery. However, once
we knew what had to be proved, a simpler argument appeared,  
yielding a more general result (\cite{LT}). Let $S$ be an arbitrary real, symmetric,
irreducible matrix with simple spectrum. Call the diagonal matrices in the closure
of $\FS$ \emph{accessible}, and their images under the BFR map the \emph{accessible vertices}.
Now let $\PS$, the \emph{spectral polytope} of $S$, be the convex closure of the 
accessible vertices of $S$.
 
\bigskip
\noindent\textbf{Theorem: } The BFR map is a diffeomorphism between the closure of the
slice $\FS$ and the spectral polytope $\PS$. 
\bigskip

Generically, the spectral polytope of a symmetric matrix with simple spectrum is
indeed the permutohedron associated to the $n!$ diagonal arrangements of the eigenvalues.
But this is not the case in general. A self-contained description of the spectral
polytope of an irreducible matrix $S$ is as follows. 
Consider the spectral decomposition,
$S = Q^T \Lambda Q$. We denote by $Q_{\{r_1,\ldots,r_k\},\{c_1,\ldots,c_k\}}$ the minor
of $Q$ consisting of the entries in the intersection of rows with indices
$r_1,\ldots, r_k$ with columns with indices in $c_1,\ldots, c_k$. For a 
permutation $pi$ in $n$ symbols, let $\Pi$ denote
the permutation matrix with entries $\Pi_{i,j} = \delta_{i,\pi(j)}$.

\bigskip
\noindent\textbf{Theorem: } A diagonal matrix $\Lambda_\pi = \Pi^T \Lambda \Pi$ is an accessible
vertex of the slice $\FS$ if and only if the minors
$Q_{\{\pi(1)\},\{1\}}, Q_{\{\pi(1),\pi(2)\},\{1,2\}},\ldots,
Q_{\{\pi(1),\ldots,\pi(n)\},\{1,\ldots,n\}}$ have nonzero determinant.
\bigskip

As an example, let $S = Q^T \Lambda Q$ be a $3 \times 3$ matrix with eigenvalues 
$\Lambda = \diag(4,2,1)$, and
so that the entry $Q_{1,1}$ equals zero, but no other minor of $Q$ has determinat
equal to zero. The slice $\FS$ and the spectral polytope have only four vertices:
diagonal matrices with entry $(1,1)$ equal to 4 do not belong to 
$\bar{\FS}$, and the vertices of the spectral polytope are $(2,1,4), (2,4,1),
(1,2,4)$ and $(1,4,2)$ --- $\bar{\FS}$ is a quadrilateral.

\bibliographystyle{amsplain}

\end{document}